\documentstyle[12pt]{article}

\newtheorem{theo}{Theorem}
\newtheorem{prop}[theo]{Proposition}
\newtheorem{coro}[theo]{Corollary}
\newtheorem{lemme}[theo]{Lemma}
\newtheorem{defi}{Definition}

\def\ptf{{X_{\tau }}}
\def\pl{1}
\def\sub{\tau = (\tau , A , \pl)}
\def\suba{\tau = (\tau , A , a)}

\def\nat{\hbox{\sf I\hskip -2pt N}}

\def\pnat{\hbox{{\scriptsize \sf I}\hskip -1pt {\scriptsize \sf N}}}

\def\card{{\rm Card}}
\def\min{{\rm min}}
\def\max{{\rm max}}
\def\inf{{\rm Inf}}

\def\mvide{{\epsilon }}

\def\der{{\cal D}}

\def\ere{{\cal R}}

\begin{document}

\begin{center}
{\LARGE\bf A generalization of Cobham's Theorem}
\end{center}

\vspace{1cm}

\begin{center}
{\large\bf Fabien Durand}
\end{center}
\begin{center}
Institut de Math\'ematiques de Luminy UPR 9016-CNRS\\
163, avenue de Luminy - Case 930 - 13288 MARSEILLE cedex~9\\
August, 28, 1997
\end{center}

\vspace{0,6cm}

\noindent
{\bf Abstract} If a non-periodic sequence $X$ is the image by a morphism of a fixed point of both a primitive substitution $\sigma$ and a primitive substitution $\tau$, then the dominant eigenvalues of the matrices of $\sigma$ and of $\tau$ are multiplicatively dependent. This is the way we propose to generalize Cobham's Theorem.

\section{Introduction}
In 1969 A. Cobham \cite{Co1} proved the following result (latter called Cobham's Theorem): {\it Let $p$ and $q$ be multiplicatively independent positive integers and $E$ a subset of $\nat$. The set $E$ is recognizable by both a $p$-automaton and a $q$-automaton if and only if $E$ is ultimately periodic.}

Later A. Cobham \cite{Co2} showed that a subset $E$ of $\nat$ is $p$-recognizable for some integer $p$ if and only if its characteristic sequence $(x_n)_{n\in \pnat}$ (i.e. $x_n = 1$ if $n$ belongs to $E$ and 0 otherwise) is $p$-substitutive (i.e. the image by a letter to letter morphism of a fixed point of a substitution of constant length $p$). There are several equivalent definitions of $p$-substitutive sequences, see for instance \cite{Co2}, \cite{CKMR} and \cite{Al}.

Hence Cobham's Theorem can be formulated as follows: {\it Let $p$ and $q$ be multiplicatively independent positive integers and $X$ be a sequence on a finite alphabet. The sequence $X$ is both $p$-substitutive and $q$-substitutive if and only if it is ultimately periodic. }

A classical result concerning matrices asserts that a square matrix with non-negative coefficients always has a real eigenvalue which is larger (not necessarily strictly) than the modulus of all other eigenvalues of $M$; moreover such an eigenvalue is a Perron number (see for instance \cite{LM}). We say that it is the dominant eigenvalue of $M$.

One can check that the dominant eigenvalue of the matrix of a substitution of constant length $p$ is $p$. To extend the notion of $p$-substitutive sequences we will say that a sequence is $\alpha$-substitutive if and only if it is the image by a letter to letter morphism of a fixed point of a substitution $\sigma$ such that $\alpha$ is the dominant eigenvalue of the matrix of $\sigma$.

Consequently a natural generalization of Cobham 's Theorem is:

\vspace{0,3cm}

{\it Let $\alpha$ and $\beta$ be two multiplicatively independent Perron numbers and $X$ a sequence on a finite alphabet. The sequence $X$ is both  $\alpha$-substitutive and $\beta$-substitutive if and only if it is ultimately periodic.}

\vspace{0,3cm}

An answer to this conjecture has been given by S. Fabre \cite{Fa1,Fa2,Fa3} in the case where $\alpha$ is a pisot number and $\beta$ a positive integer. Recently, using the formalism of the first order logic A. B\`es \cite{Be} and I. Fagnot \cite{Fag} obtained a partial answer in the case where $\alpha$ and $\beta$ are pisot numbers.

In this paper we give a positive answer to this conjecture in the case where the substitutions are primitive without any assumption concerning the eigenvalues $\alpha$ and $\beta$:

\vspace{0,3cm}

{\it Let $\alpha$ and $\beta$ be two multiplicatively independent Perron numbers and $X$ a sequence on a finite alphabet. If a sequence $X$ is both  $\alpha$-substitutive and $\beta$-substitutive then $X$ is periodic.}

\vspace{0,3cm}

There are other ways to generalize Cobham's Theorem, some of them can be found in \cite{Be}, \cite{Fa2}, \cite{Fag}, \cite{MV2} and \cite{Se}. Related works can be found in \cite{BHMV}, \cite{Ha} or \cite{MV1}.

\vspace{0,3cm}

Section 2 of this paper contains the basic definitions we need. In Section 3 we define the main notion of this paper, the {\it return word}, which was first introduced in \cite{Du}. We review some properties of return words obtained in \cite{Du}. Section 4 is an intermediate step to prove the main result where we establish helpful morphism relations. Section 5 is split into two subsections. In the first one we prove a result (Theorem \ref{propvalpr}) stronger than our main theorem though only valid for fixed points: 

\vspace{0,3cm}

\noindent {\it If two primitive substitutions have the same non-periodic fixed point, then they have some powers which have the same eigenvalues, except perhaps 0 and the roots of the unity.}

\vspace{0,3cm}

 In the second one we prove the main theorem. An example is used to show that for substitutive primitive sequences we cannot have a better result. The aim of Section 6 is to show that there is more than only relations between eigenvalues. Primitive substitutions sharing a same fixed point have some powers which coincide on some sets of return words.

\section{Definition and terminology}

\subsection{Words and sequences}

We call {\it alphabet} a finite set of elements called {\it letters}. Let $A$ be an alphabet, a {\it word} on $A$ is an element of the free monoid on $A$, denoted by $A^*$, i.e. a finite (possibly empty) sequence of letters. Let $x = x_0x_1 \cdots x_{n-1}$ be a word, its {\it length} is $n$ and is denoted by $|x|$. The {\it empty-word} is denoted by $\mvide$, $|\mvide| = 0$. The set of non-empty words on  $A$ is denoted by $A^+$. The elements of $A^{\pnat}$ are called {\it sequences}. If $X=X_0X_1\cdots$ is a sequence (with $X_i\in A$, $i\in \nat$), and $l$, $k$ are two non-negative integers, with $l\geq k$, we denote the word $X_k X_{k+1}\cdots X_{l}$ by $X_{[k,l]}$ and we say that $X_{[k,l]}$ is a {\it factor} of $X$. If $k = 0$, we say that $X_{[0,l]}$ is a {\it prefix} of $X$ and we write $X_{[0,l]} \prec X$. The set of factors of length $n$ of $X$ is written $L_n(X)$, and the set of factors of $X$, or {\it language} of $X$, is represented by $L(X)$. If $u$ is a factor of $X$, we will call {\it
occurrence} of $u$ in $X$ every integer $i$ such that $X_{[i,i + |u| - 1]}= u$. When $X$ is a word, we use the same terminology with the
similar definitions. Let $u$ and $v$ be two words, we denote by $L_u(v)$ the number of occurrences of $u$ in $v$. A word $u$ is a {\it suffix} of the word $v$ if $v=xu$ for some $x$ belonging to $A^*$.

The sequence $X$ is {\it ultimately periodic} if there
exist a word $u$ and a non-empty word $v$ such that $X=uv^{\omega}$, where $v^{\omega}$
is the infinite concatenation of the word $v$. Otherwise we say that $X$ is {\it non-periodic}. It is {\it periodic} if $u$ is the empty-word. 

A sequence $X$ is {\it uniformly recurrent} if for each factor $u$ the greatest difference of two successive occurrences of $u$ is bounded.

\subsection{Morphisms and matrices}

Let $A$, $B$ and $C$ be three alphabets. A {\it morphism} $\tau$ is a map from $A$ to $B^*$. Such a map induces by concatenation a map from $A^*$ to $B^*$. If $\tau (A)$ is included in $B^+$, it induces a map from $A^{\pnat}$ to $B^{\pnat}$. All these maps are written $\tau$ also.

To a morphism $\tau$, from $A$ to $B^*$, is naturally associated the matrix $M_{\tau} = (m_{i,j})_{i\in  B  , j \in  A  }$ where $m_{i,j}$ is the number of occurrences of $i$ in the word $\tau(j)$. To the composition of morphisms corresponds the multiplication of matrices. For example, let $\tau_1 : B \rightarrow C^*$, $\tau_2 : A \rightarrow B^*$ and $\tau_3 : A \rightarrow C^*$ be three morphisms such that $\tau_1 \tau_2 = \tau_3$, then we have the following equality: $M_{\tau_1} M_{\tau_2} = M_{\tau_3}$. In particular if $\tau $ is a morphism from $A$ to $A^{*}$ we have $M_{\tau^n} = M_{\tau}^n$.

A non-negative square matrix $M$ always has a non-negative eigenvalue $r$ such that the modulus of all its other eigenvalues do not exceed $r$. We call it the {\it dominant eigenvalue of $M$} (see for instance \cite{Ga}). A  square matrix is called {\it primitive} if it has a power with positive coefficients. A morphism from $A$ to $A^*$ is called {\it primitive} if its associated matrix is primitive. In this case the dominant eigenvalue is a simple root of the characteristic polynomial, and is strictly larger than the modulus of all other eigenvalues. This is Perron's Theorem (\cite{Ga}, p. 53).

\subsection{Substitutions and substitutive sequences}

\begin{defi} 
A {\rm substitution} is a triple $\tau = (\tau , A ,
a)$, where $A $ is an alphabet, $\tau $ is a morphism from $A $ to $A^+$ and $a$ is a letter of $A$ such that the first letter of $\tau (a )$ is $a$.
\end{defi}

Let $\tau = (\tau , A ,a )$ be a substitution. There exists a unique sequence $X=(x_n)_{n\in \pnat}$ of $A^{\pnat}$ such that $x_0 = a$ and $\tau (X) = X$ (for more details we refer the reader to \cite{Qu}). We will say that $X$ is the {\it fixed point} of $\tau$ and we will denote it by $X_{\tau}$.

In this article we only consider {\it primitive substitutions}, i.e. substitutions with primitive associated matrices. If $\tau = (\tau , A ,
a)$ is a primitive substitution it is not difficult to see that its fixed point is uniformly recurrent (see \cite{Qu}).

Let $A$ and $B$ be two alphabets, we say that a morphism $\sigma$ from $A$ to $B^*$ is a {\it letter to letter morphism} when $\sigma (A)$ is a subset of $B$. A sequence $Y$ is {\it substitutive} if there exist a primitive substitution $\tau$ and a letter to letter morphism $\phi$ such that $Y = \phi (\ptf)$. We will also say that $Y$ arises from $\tau$. It is {\it $\alpha$-substitutive} if $\alpha$ is the dominant eigenvalue of $\tau$. We can remark that each substitutive sequence is uniformly recurrent. In particular, it is periodic whenever it is ultimately periodic.

\vspace{0,3cm}

From the proof of Proposition 9 in \cite{Du} we deduce the following proposition:
\begin{prop}
Let $A$ and $B$ be two alphabets, $X$ be a $\alpha$-substitutive sequence on $A$ and $\varphi : A\rightarrow B^+$ be a morphism. There exists a positive integer $k$ such that the sequence $\varphi (X)$ is $\alpha^k$-substitutive.
\end{prop}
This proposition allows us to consider only letter to letter morphisms without loss of generality.

\section{Return words}

In this section we define the main notion used in this paper, the {\it return words}. It was introduced in \cite{Du} where we stated and proved some of its properties we recall here. We will use them very often in the sequel.

\subsection{Definition}

Let $X$ be a uniformly recurrent sequence on the alphabet $A$ and $u$ a non-empty prefix of $X$. We call {\it return word on u} every factor
$X_{[i,j-1]}$, where $i$ and $j$ are two successive occurrences of $u$ in $X$. For example let 
$$
X = ababcababbbabababcababbbababaccababacc \cdots
$$
be a sequence. The words $ababc$, $ababbb$, $ab$, $ababacc$ are return words on $abab$ of $X$.

The reader can check that a word $v$ is a return word on $u$ of $X$ if and only if $vu$ belongs to $L(X)$, $u$ is a prefix of $vu$ and $u$ has exactly two occurrences in $vu$. For the details we refer the reader to \cite{Du}. The set of return words on $u$ is finite, because $X$ is uniformly recurrent, and is denoted by ${\ere}_{X,u}$. The sequence $X$ can be written naturally as a concatenation
$$
X = m_0m_1m_2 \cdots \:\: , m_i \in \ere_{X,u}, \:\: i\in \nat,
$$
of return words on $u$, and this decomposition is unique. We enumerate the elements of $\ere_{X,u}$ in the order of their first appearance in $(m_n)_{n\in\pnat}$. This defines a bijective map
$$
\Theta_{X,u} : R_{X,u}  \rightarrow {\ere}_{X,u} \subset A^{*}
$$
where $R_{X,u} = \{ 1, \cdots , \card({\ere}_{X,u})\}$. 

The map $\Theta_{X,u}$ defines a morphism and the set $\Theta_{X,u} (R_{X,u} ^*)$ consists of all concatenations of return words on $u$.
We denote by $\der_u(X)$ the unique sequence on the alphabet $R_{X,u}$ characterized by
$$
\Theta_{X,u} (\der_u(X) ) = X.
$$
We call it the {\it derived sequence of $X$ on $u$}. It is clearly uniformly recurrent. We remark that
$$
L(X) \bigcap \Theta_{X,u} (R_{X,u} ^*) = \Theta_{X,u} (L(\der_u(X))).
$$
When it does not create confusion we will forget the ``$X$'' in the symbols $\Theta_{X,u}$, $\ere_{X,u}$ and $R_{X,u}$.

\subsection{Some properties of return words}

The following proposition points out the basic properties of return words which are of constant use throughout the paper.

\begin{prop}
\label{derder}
{\rm (Proposition 6 in \cite{Du})} Let $X$ be a uniformly recurrent sequence and $u$ a non-empty prefix of $X$.
\begin{enumerate}
\item The set $\ere_{X,u}$ is a code, i.e. $\Theta_{X,u} : R_{X,u}^* \rightarrow \Theta_{X,u} (R_{X,u} ^*)$ is one to one. 
\item If $u$ and $v$ are two prefixes of $X$ such that $u$ is a prefix of $ v$ then each return word on $v$ belongs to $\Theta_{X,u} (R_{X,u} ^*)$, i.e. it is a concatenation of return words on $u$.
\item Let $v$ be a non-empty prefix of $\der_u (X)$ and $w = \Theta_{X,u}(v)u$. Then 
\begin{itemize}
\item $w$ is a  prefix $X$,
\item $\Theta_{X,u}\Theta_{\der_u(X),v} = \Theta_{X,w}$ and
\item $\der_v (\der_u (X)) = \der_w(X)$. 
\end{itemize}
\end{enumerate}
\end{prop}

\begin{lemme}
\label{lemtech}
{\rm (Lemma 10 in \cite{Du})} Let $X$ be a non-periodic uniformly recurrent sequence, then
$$
m_n = \inf \{ | v | ; v\in \ere_{X,X_{[0,n]}} \} \rightarrow +\infty \;\; {\rm when} \;\; n\rightarrow +\infty.$$
\end{lemme}

\subsection{Substitutive sequences and return words}

When we apply return words to primitive substitutions we obtain some useful results. The following proposition states that each derived sequence of a fixed point of a primitive substitution is a fixed point of a primitive substitution too.

\begin{prop}
\label{fpderived}
{\rm (Proposition 19 in \cite{Du})} Let $\suba $ be a primitive substitution and $u$ be
a non-empty prefix of $\ptf $. The derived sequence $\der_u (\ptf) $ is the
fixed point of a primitive substitution $\tau_{u} = (\tau_u ,R_{u}, 1)$ where $\tau_u$ satisfies
$$
\Theta_u \tau_u = \tau \Theta_u.
$$  
\end{prop}

The map $\Theta_u$ being one to one the previous equality completely characterized $\tau_u$. Such a substitution is called {\it return substitution} ({\it on $u$}). Moreover we can remark that $(\tau^l)_u = (\tau_u)^l$.

The two following theorems were established in \cite{Du} to obtain a characterization of substitutive sequences: A non-periodic uniformly recurrent sequence $Y$ is substitutive if and only if the set of its derived sequences is finite.
\begin{theo}
\label{encad}
{\rm (Theorem 18 in \cite{Du})} Let $Y$ be a non-periodic substitutive sequence. There exist three positive constants $H_1$, $H_2$ and $H_3$ such that: For all non-empty
prefixes $u$ of $Y$, 
\begin{enumerate}
\item for all words $v$ belonging to $\ere_{Y,u}$, $H_1| u| \leq | v| \leq H_2| u|$, and 
\item $\card (\ere_{Y,u})\leq H_3$.
\end{enumerate}
\end{theo}

\begin{theo}
\label{subretfin}
{\rm (Theorem 20 in \cite{Du})} Let $\suba$ be a primitive substitution. The set of the return substitutions of $\tau$ is finite.
\end{theo}

\section{Eigenvalues and return words}

We establish some morphism relations between the substitutions and their return substitutions, then we find their common eigenvalues.

In this section $\suba$ will be a primitive substitution and $u$, $v$ two prefixes of $\ptf$ such that $|u| < |v|$. We recall that we have
\begin{eqnarray}
\label{eq1}
\Theta_u \tau_u = \tau \Theta_u \:\: {\rm and} \:\: \Theta_v \tau_v = \tau 
\Theta_v .
\end{eqnarray}
The word $u$ is a prefix of $v$, hence a return word on $v$ is a concatenation of return words on $u$. This allows us to define the morphism $\lambda $, from $R_v$ to $R_u^{+}$, by $\Theta_u \lambda = \Theta_v$.
Thus we obtain the relation
$$
\tau_u \lambda = \lambda \tau_v.
$$

Let $k$ be an integer such that $|v| < |\tau^k (u)|$. The image by $\tau^k$ of a return word on $u$ is a concatenation of return words on $v$. We define a new morphism $\kappa$, from $R_u$ to $R_v^{+}$,
by $\Theta_v \kappa = \tau^k \Theta_u$. We deduce the following morphism relations:
\[
\begin{array}{cccc}
\tau_v \kappa  & = & \kappa \tau_u, &\\
\kappa \lambda & = & \tau_v^k & \:\: {\rm and}\\
\lambda \kappa & = & \tau_u^k. &
\end{array}
\]
Consequently we have the following proposition:

\begin{prop}
\label{propprec}
Let $\suba$ be a primitive substitution and $u$, $v$ be two prefixes of $\ptf$ such that $|u| < |v|$. Then there exist an integer $k\geq 1$ and two morphisms $\lambda : R_v \rightarrow R_u^{+}$ and $\kappa : R_u \rightarrow R_v^{+}$ such that
$$
\tau_v \kappa = \kappa \tau_u ,\:\: \tau_u \lambda = \lambda \tau_v,  \:\:
\kappa \lambda = \tau_v^k \:\: {\rm and} \:\: \lambda \kappa = \tau_u^k .
$$
\end{prop}

\begin{coro}
All the return substitutions of a primitive substitution have all the same non-zero eigenvalues.
\end{coro}
Proof: This is a straightforward consequence of Proposition \ref{propprec}. The details are left to the reader. \hfill $\Box$

\vspace{0,3cm}

By primitivity, there exists an integer $n_0$ such that for all $n\geq n_0$ all images by $\tau^{n}$ of letters have at least two occurrences of $u$. Let $l$ be an integer larger than $n_0$ and $K_l$ be the matrix defined by
$$
K_l = (L_{\Theta_u (c)u}(\tau^l (b)u) )_{c\in R_u,b \in A};
$$ 
we recall that $L_{\Theta_u (c)u}(\tau^l (b)u)$ is the number of occurrences of $\Theta_u(c)u$ in $\tau^l (b)u$. 

Let $b$ be an element of $A$. We set $\tau^l(b)u = a_0 a_1 \cdots a_k a_{k+1} \cdots a_{k+|u|}$. Let $i$ be the first occurrence of $u$ in $\tau^l(b)$ and $j$ the greatest occurrence of $u$ in $\tau^l(b)u$ such that $a_i \cdots  a_{j-1}$ is a concatenation of elements of $\ere_{X_{\tau},u}$. We set $x=a_0 a_1\cdots a_{i-1}$, $y=a_j a_{j+1}\cdots a_{k+|u|}$ and $w = a_i a_{i+1}\cdots a_{j-1}$. The word $w$ is a concatenation of return words on $u$ and $L_{\Theta_u (d)u}(wu) = L_{\Theta_u (d)u}(\tau^l (b)u)$. We remark that the length of $x$ is less than $H_2 |u|$ and that the length of $y$ is less than  $(H_2 +2)|u|$ where $H_2$ is the constant given by Theorem \ref{encad}.

Let $c$ be a letter of $A$,
$$
L_c (\tau^l (b)) = L_c (x) + \sum_{d\in R_u} L_c (\Theta_u (d)) L_{\Theta_u (d)u}(wu)  + L_c (y) 
$$
$$ = L_c (x) + \sum_{d\in R_u} L_c (\Theta_u (d)) L_{\Theta_u (d)u}(\tau^l (b)u)  + L_c (y).
$$

We observe that the number of occurrences of $c$ in both $x$ and $y$ is less $(H_2 +2)|u|$. Then we have
\begin{eqnarray}
\label{eq6}
M_{\tau^l} = M_{\tau}^l = M_{\Theta_u} K_l + Q_l,
\end{eqnarray}
where $Q_l$ is a non-negative integral matrix whose coefficients are less than $(H_2 + 2) |u|$. For this reason the set $\{ Q_l ; l\in \nat \}$ is finite.

Let $b$ and $c$ be two elements of $R_u$. We set 
$$
p_{c,b} =
L_{\Theta_u (c)u}(\tau^l (\Theta_u (b))u) - \sum_{d\in A} L_{\Theta_u (c)u} ( \tau^l (d)u) L_d(\Theta_u (b)) \ .
$$
We can bound $p_{c,b}$ independently of $l$. Let $xy$ be a word of length 2 occurring in $\Theta_u (b)$. We set $\tau^l (xy) = v_0 \cdots v_{k-1} v_{k} \cdots v_n$ where $k = |\tau^l (x)|$. Let $j$ be the greatest occurrence of $u$ in $\tau^l (xy)$ less or equal to $k-1$ and $i$ be the smallest occurrence of $u$ in $\tau^l (xy)$ larger or equal to $k$. We have
$$
|L_{\Theta_u (c)u} (v_i \cdots v_{j-1}u) - (L_{\Theta_u (c)u} (v_i\cdots v_{k-1} u) + L_{\Theta_u (c)u} (v_k\cdots v_{j-1} u )) | \leq 
$$
$$
\frac{2(H_2 +1)|u|}{H_1|u|} = \frac{2(H_2 +1)}{H_1}.
$$
Where $H_1$ is the constant given by Theorem \ref{encad}. Hence
$$
p_{c,b} \leq \frac{2(H_2 +1)}{H_1} |\Theta_u (b)| \leq \frac{2(H_2 +1)H_2 |u|}{H_1}.
$$
Moreover we have
$$
L_c (\tau^l_u (b)) = L_{\Theta_u (c)u} ( \tau^l (\Theta_u (b))u) = \sum_{d\in A} L_{\Theta_u (c)u} ( \tau^l (d)u) L_d(\Theta_u (b)) + p_{c,b} \ .
$$
Consequently we have 
\begin{eqnarray}
\label{eq7}
M_{\tau_u^l} = M_{\tau_u}^l = K_l M_{\Theta_u} + P_l \ .
\end{eqnarray}
where $P_l$ is an integral matrix; the absolute values of its coefficients are less than $2(H_2 + 1)H_2|u|/H_1$.
Therefore the set $\{ P_l ; l\in \nat \}$ is finite.

\begin{prop}
\label{relvalpr}
Let $\tau$ be a primitive substitution and $u$ a prefix of $\ptf$. The substitutions $\tau $ and $\tau_u$ have the same eigenvalues, except perhaps $0$ and roots of the unity.
\end{prop}
Proof: Let $u$ be a prefix of $\ptf$ and $\alpha$ a non-zero eigenvalue of $M_{\tau}$ which is not a root of the unity. There exists a vector $x \not = 0$ such that $(^T x)M_{\tau} = \alpha (^T x)$. According to the relation (\ref{eq1}), we have 
$$
(^T
x)M_{\Theta_u} M_{\tau_u} = \alpha (^T x)M_{\Theta_u}.
$$
We have to prove that $(^T x)M_{\Theta_u}$ is different from zero.

Suppose it is false. From equality (\ref{eq6}) it follows that $\alpha^l (^T x) = (^T x)Q_l$ for all integers $l$ larger than $n_0$. But the set $\{ Q_l, l\in \nat \}$ is finite. Hence there exist two distinct integers $l_1$ and $l_2$, larger than $n_0$, such that $Q_{l_1} = Q_{l_2}$. And finally we have $\alpha = 0$ or $\alpha^{l_1-l_2} = 1$, which contradicts our assumption on $\alpha$.

In the same way, it follows from equality (\ref{eq7}) that if $\mu$ is a non-zero eigenvalue of $M_{\tau_u}$ which is not a root of the unity, then $\mu$ is an eigenvalue of $M_{\tau}$. This completes the proof.
\hfill $\Box$

\vspace{0,3cm}

It is easy to check that if $\tau : \{0,1 \} \rightarrow \{0,1 \}^+$ is the Fibonacci substitution, i.e. $\tau (0) = 01$ and $\tau(1) = 0$, then we have $\tau = \tau_{01}$. Hence $\tau$ and $\tau_{01}$ have the same eigenvalues. On the other hand the set of eigenvalues of the Morse substitution, $\sigma(0)=01$ and $\sigma(1)=10$, is $\{ 0,2 \}$ and the eigenvalues of $\sigma_{011}$ are $0,0,-1$ and $2$.

\section{A generalization of Cobham's Theorem}

The proof of the generalization we announced requires several steps. In Proposition \ref{stepone} and Proposition \ref{corotech} we work under special assumptions. Proposition \ref{corotech} shows why these assumptions are relevant for our purpose and Lemma \ref{lemmetech} proves that it is always possible to work under these assumptions. This leads to a stronger theorem than the generalization though only valid for fixed points as announced in the introduction.

For convenience in the sequel we will use alphabets $\{ 1,2,\cdots ,k \}$.

\subsection{Some technical results}

\begin{prop}
\label{stepone}
Let $\sub$ be a primitive substitution and $u$ be a prefix of $\ptf$ such that:
\begin{enumerate}
\item 
\label{hypzero}
For all letters $b$ of $A$, $\tau(b)$ begins by $\pl$,
\item
\label{hypdeux}
The substitutions $\tau$ and $\tau_u$ are defined on the same alphabet and are identical, 
\item 
\label{hyptrois} 
The fixed point of $\tau$ is non-periodic,
\item 
\label{hypquatre}
For all letters $b$ and $c$ of $A$, $b$ has at least one occurrence in $\Theta_u(c)$.
\end{enumerate}
Let $J$ be an infinite set of positive integers. Then there exist an infinite subset $I$ of $J$, a strictly increasing sequence of positive integers $(l_p)_{p\in I}$ and a morphism $\gamma : A \rightarrow A^+$ such that for all $p$ in $I$.
$$
\Theta_u^{l_p} \gamma = \gamma \Theta_u^{l_p} = \tau^p .
$$
\end{prop}
Proof: Hypothesis \ref{hypdeux} says that $A=R_u$. It is easy to check that the morphism $\Theta_u : A\rightarrow A^*$ defines a substitution $\Theta_u = (\Theta_u ,A , \pl)$. We put $\Theta = \Theta_u$. Hypothesis \ref{hypquatre} implies that this substitution is primitive. 

As the substitutions $\tau_u$ and $\tau$ are identical (hypothesis \ref{hypdeux}), they have the same fixed point $X_{\tau}$ and we have seen that the fixed point of $\tau_u$ is $\der_u (\ptf)$ (Proposition \ref{fpderived}), hence $\der_u (\ptf) = \ptf$. Consequently, we have $\ptf = \Theta(\ptf)$, i.e. $X_{\tau}$ is the fixed point of $\Theta = (\Theta ,A , 1)$. 

Moreover we can remark that $\tau \Theta = \Theta \tau$.

The word $u$ is a prefix of $\der_u(\ptf)$, hence we can consider the sequences $(\der_u^n(\ptf))_{n\geq 1}$ defined by
\begin{center}
$\der_u^1 (\ptf) = \der_u(\ptf)$ and $ \der_u^{n+1} (\ptf) = \der_u (\der_u^n (\ptf))$ for all $n\geq 1$.
\end{center}

Let us prove by induction that for all $n\geq 1$ we have:
\begin{enumerate}
\item[i)] 
$\der_u^n (\ptf) = \der_{w_n} (\ptf) = \ptf$, with $w_n = \Theta^{n-1}(u) \cdots \Theta (u) u$,
\item[ii)]
$\Theta^n = \Theta_{w_n}$ and
\item[iii)] 
$\tau = \tau_{w_n}$.
\end{enumerate}

For $n=1$ it suffices to remark that $w_1 = u$. 

Now suppose that points i), ii) and iii) are satisfied for some positive integer $n$. We have 
$$
\der_u^{n+1} (\ptf) = \der_u (\der_u^{n}(\ptf)) = \der_u (\ptf) = X_{\tau}
$$
and Proposition \ref{derder} implies that:
\begin{itemize}
\item
$\der_u^{n+1} (\ptf) = \der_u (\der_{w_n}(\ptf)) = \der_w (\ptf)$ and
\item
$\Theta_{w} = \Theta_{w_n} \Theta_{\der_{w_n} (X_{\tau}),u} = \Theta^n \Theta_{X_{\tau} , u} = \Theta^{n+1}$
\end{itemize}
where 
$$
w = \Theta_{w_n} (u) w_n = \Theta^n (u) \Theta^{n-1} (u) \cdots \Theta (u) u = w_{n+1}.
$$ 
Hence points i), ii) are satisfied for $n+1$.

The substitution $\tau_{w_{n+1}}$ is the return substitution on $u$ of $\tau_{w_n}$ consequently $\Theta \tau_{w_{n+1}} = \tau_{w_n} \Theta$; that is to say $\Theta \tau_{w_{n+1}} = \tau \Theta$. But $\Theta \tau_u = \tau \Theta$ and the map $\Theta$ is one to one hence $\tau_{w_{n+1}} = \tau_u = \tau$. This completes the proof by induction of points i), ii) and iii).

We denote the dominant eigenvalues of $M_{\tau}$ and $M_{\Theta}$  respectively by $\alpha$ and $\beta$. We recall (see for instance \cite{Qu}) that there exists a positive number $r$ such that for all $b$ in $A$ and all $k$ in $\nat$
$$
\frac{1}{r} \alpha^k \leq | \tau^k (b) | \leq r\alpha^k \:\: 
{\rm and} \:\: \frac{1}{r} \beta^k \leq | \Theta^k (b) | \leq 
r\beta^k .
$$
From this we deduce that there exists two constants $c_1$ and $c_2$ such that for all positive integers $n$ 
$$
c_1 \beta^n \leq | w_n | = |\Theta^{n-1}(u) \cdots \Theta (u) u| \leq c_2 \beta^n.
$$
From hypothesis \ref{hypzero} it follows that there exists an integer $k_0$ such that $u$ is a prefix of all images of letters by $\tau^{k_0}$. For every integer $k$, larger than $k_0$, we define $l_k$ to be the greatest integer $n$ such that $w_n$ is a prefix of $\tau^{k-1}(1)$. For all positive integers we have 
$$
c_1 \beta^{l_k} \leq | w_{l_k} | \leq | \tau^{k-1} (1)| \leq | w_{l_{k}+1} | \leq c_2 \beta^{l_k +1}.
$$
Thus we obtain 
$$
| \tau^{k-1} (1)| \leq \frac{\beta c_2}{c_1} c_1 \beta^{l_k} \leq \frac{\beta c_2}{c_1} | w_{l_k} |.
$$

Let $k$ be an integer larger than $k_0$. For all letters $b$ of $A$ the word $w_{l_k}$ is a prefix of $\tau^k (b)$. Hence all images by $\tau^k$ of words are concatenations of return words on $w_{l_k}$. This remark allows us to define the morphism $\gamma_k$, from $A$ to $A^+$, by $\Theta_{w_{l_k}} \gamma_k = \tau^k$. We have:
$$
\Theta_{w_{l_k}} \gamma_k \Theta_{w_{l_k}} = \tau^k \Theta^{l_k} = \Theta^{l_k} \tau^k.
$$
The map $\Theta_{w_{l_k}} = \Theta^{l_k}$ is one to one hence $\gamma_k \Theta_{w_{l_k}} = \tau^k$ and finally
\begin{eqnarray}
\label{rel1}
\Theta_{w_{l_k}} \gamma_k = \gamma_k \Theta_{w_{l_k}} .
\end{eqnarray}
Moreover for all $b$ in $A$
$$
|\gamma_k (b)| = L_{w_{l_k}} (\tau^k (b)w_{l_k}) - 1
                      \leq \frac{|\tau^k (b)| + |w_{l_k}| }{H_1|w_{l_k}|}                      
                      \leq \frac{2\beta c_2 | \tau^k (b)|}{H_1 c_1 | \tau^{k-1} (1) | } 
                      \leq \frac{2c_2 r^2 \alpha \beta}{H_1 c_1 },
$$
where $H_1$ is the constant of Theorem \ref{encad}. We have proved that the length of the images by $\gamma_k$ of letters are bounded independently of $k$. Hence the set $\{ \gamma_k ; k\geq k_0 \}$ is finite. Thus there exists an infinite set $I$, included in $J$, such that $\gamma_p = \gamma_q$ for all elements $p$ and $q$ of $I$.

Let $p$ be an element of $I$, we write $\gamma = \gamma_p$. Equality (\ref{rel1}) gives $\Theta^{l_p} \gamma = \gamma \Theta^{l_p} = \tau^p$. From this last equality it follows that the sequence $(l_p)_{p\in I}$ is strictly increasing.\hfill $\Box$

\vspace{0,3cm}

\begin{prop}
\label{corotech}
Let $\sub$ be a primitive substitution and $u$ be a prefix of $\ptf$ satisfying the hypothesis of Proposition \ref{stepone}. Let $J$ be an infinite set of positive integers. Then there exist an infinite subset $I$ of $J$ and a strictly increasing sequence of positive integers $(l_p)_{p\in I}$ such that, for all distinct integers $p$ and $q$ belonging to $I$, $p<q$, we have:
\begin{itemize}
\item 
The non-zero eigenvalues of $\tau^{q-p}$ are eigenvalues of $\Theta_u^{l_q-l_p}$,
\item 
$\Theta_u$ is a primitive substitution and the non-zero eigenvalues of $\Theta_u^{l_q-l_p}$, except perhaps roots of the unity, are eigenvalues of $\tau^{q-p}$.
\end{itemize}
\end{prop}
Proof: We set $\Theta = \Theta_u$. As in the proof of Proposition \ref{stepone} we can remark that $\Theta$ defines a primitive substitution $\Theta = (\Theta , A , 1)$. There exist an infinite set $I$ of integers, a strictly increasing sequence of integers $(l_p)_{p\in I}$ and a morphism $\gamma$ from $A$ to $A^{+}$ such that for all $p$ in $I$
$$
\Theta^{l_p} \gamma = \gamma \Theta^{l_p} = \tau^p.
$$
Let $p<q$ be two elements of $I$. From the previous equalities and from the fact that $\Theta \tau = \tau \Theta$ (because $\tau = \tau_u$) we obtain
\begin{eqnarray}
\label{relgamma}
\tau^q = \Theta^{l_q - l_p} \Theta^{l_p} \gamma = \Theta^{l_q - l_p} \tau^p = \tau^p \Theta^{l_q - l_p}.
\end{eqnarray}

Let $\alpha $ be a non-zero eigenvalue of $M_{\tau^{q-p}} = M_{\tau}^{q-p}$ and $x$ one of its eigenvectors. Equality (\ref{relgamma}) implies 
$$
M_{\Theta}^{l_q-l_p} M_{\tau}^p x  = M_{\tau}^{q}x = \alpha M_{\tau}^p x  .
$$ 
But the vector $M_{\tau}^{(q-p)p} x$ is different from zero and therefore so is $M_{\tau}^p x$. Thus $\alpha$ is an eigenvalue of $\Theta^{l_q - l_p}$. This completes the first part of the proof.

\vskip 0,3cm

It remains to prove the second part. Let $k_0$ be an integer such that $\min_{b\in A} |\Theta^{k_0} (b)| > 2\max_{b\in A} |\tau^p (b)|$. Let $k$ be a positive integer such that $k(l_q - l_p) > k_0$. For all letters $b$ of $A$ we can write $\Theta^{k(l_q-l_p)}(b) = u\tau^p (w) v$, where $u$ and $v$ are respectively a suffix and a prefix of an element of $\tau^p (A)$, and $w$ is a non-empty word of $L(\ptf)$. Hence there exist two integral matrices with non-negative coefficients, $S_k$ and $P_k$, such that
\begin{eqnarray}
\label{reltheta}
M_{\Theta}^{k(l_q - l_p)} = M_{\tau}^p S_k + P_k,
\end{eqnarray}
where the coefficients of $P_k$ are less than $2\max_{b\in A} |\tau^p (b)|$. The set $\{ P_j ; j > k(l_q - l_p) \}$ is finite because the coefficients are bounded independently of $k$. Consequently there exist two integers $k_1$ and $k_2$ such that $P_{k_1} = P_{k_2}$.

Let $\alpha$ be a non-zero eigenvalue of $ M_{\Theta}^{l_q-l_p}$ which is not a root of the unity, and $x$ one of its associated left eigenvectors. Equality (\ref{relgamma}) leads to
$$
x M_{\tau}^{p} M_{\tau}^{q-p} = x M_{\tau}^{q} = x M_{\Theta}^{l_q - l_p} M_{\tau}^{p} = \alpha x M_{\tau}^p .
$$ 
But $P_{k_1} = P_{k_2}$, so $x M_{\tau}^p  $ is different from zero. Otherwise, by equality (\ref{reltheta}) we would have $\alpha^{k_1} x = \alpha^{k_2 } x$. Which contradicts our hypothesis on $\alpha$. \hfill $\Box$

\vspace{0,3cm}

\begin{lemme}
\label{lemmetech}
Let $(\tau, A,1)$ and $\sigma = (\sigma , A , \pl)$ be two primitive substitutions having the same non-periodic fixed point $X$. There exist an integer $l$, a prefix $v$ of $X$, and an arbitrarily long prefix $u$ of $\der_v(X)$ such that the word $u$ and the substitution $\tau_v^l$ , and $u$ and the substitution $\sigma_v^l$, both satisfy the hypothesis of Proposition \ref{stepone}.
\end{lemme}
Proof: Let $( X^{(n)} )_{ n\in \pnat }$ be the sequences defined by: $X^{(0)} = X$ and $X^{(n+1)}$ is the derived sequence of $X^{(n)}$ on the prefix $\pl$ (the first letter of $X^{(n)}$). For all integers $n$ we call $A^{(n)}$ the alphabet of $X^{(n)}$. Let $(u^{(n)})_{n\geq 1}$ be the sequence of words defined by:
$$ 
u^{(1)} = 1 \ \ {\rm and} \ \ u^{(n+1)} = \Theta_{X , u^{(n)}} (1) u^{(n)}.
$$
According to Proposition \ref{derder}, for all integers $n$ larger than 1, the word $u^{(n)}$ is a prefix of $X_{\tau}$ such that
\begin{center}
$\Theta_{X , 1} \Theta_{X^{(1)} , 1} \cdots \Theta_{X^{(n-1)} , 1} = \Theta_{X , u^{(n)}} $ and $X^{(n)} = \der_{u^{(n)}}(X)$.
\end{center}
The sets $\{ \tau_u , u\prec X \}$ and $\{ \sigma_u , u\prec X \}$ are finite by Theorem \ref{subretfin}. Hence, there exists an infinite set $I$ of positive integers such that for all integers $p$ and $q$ of $I$, we have $\tau_{u^{(p)}} = \tau_{u^{(q)}}$ and $\sigma_{u^{(p)}} = \sigma_{u^{(q)}}$. We remark that the fixed point of the substitutions $\tau_{u^{(p)}} = (\tau_{u^{(p)}}, A^{(p)} , 1)$ and $\sigma_{u^{(p)}} = (\sigma_{u^{(p)}}, A^{(p)} , 1)$ is $X^{(p)} = \der_{u^{(p)}} (X)$.

Let $p$ and $q$ be two elements of $I$ with $p<q$. By definition of $(X^{(n)})_{n\in\pnat}$ we have 
$$
X^{(q)} = \underbrace{\der_1 \cdots \der_1}_{(q-p)\ \ {\rm times}}(X^{(p)})
$$ 
Hence (Proposition \ref{derder}) there exists a prefix $u$ of $X^{(p)}$ such that 
\begin{itemize}
\item 
$\der_u (X^{(p)}) = X^{(q)}$,
\item 
$ \Theta_{X^{(p)} , 1} \Theta_{X^{(p+1)} , 1} \cdots \Theta_{X^{(q-1)} , 1} = \Theta_{X^{(p)} , u}$ \ ,
\item 
$\Theta_{X^{(p)},u} \tau_{u^{(q)}} = \tau_{u^{(p)}} \Theta_{X^{(p)},u}$ and $\Theta_{X^{(p)},u} \sigma_{u^{(q)}} = \sigma_{u^{(p)}} \Theta_{X^{(p)},u}$.
\end{itemize}
From the last equalities it is clear that $(\tau_{u^{(p)}})_u =  \tau_{u^{(q)}}$ and $(\sigma_{u^{(p)}})_u =  \sigma_{u^{(q)}}$; where $(\tau_{u^{(p)}})_u$ and $(\sigma_{u^{(p)}})_u$ are respectively the return substitutions on $u$ of $\tau_{u^{(p)}}$ and $\sigma_{u^{(p)}}$.

From the definition of $(u^{(n)})_{n\geq 1}$ we deduce that the sequence $(|u^{(n)}|)_{n\geq 1}$ is strictly increasing. Thus it follows from Lemma \ref{lemtech} that 
$$
\lim_{j\rightarrow + \infty} \min \{ |\Theta_{X^{(0)},u^{(j)}} (b)| = | \Theta_{X^{(0)} , 1} \Theta_{X^{(1)} , 1} \cdots \Theta_{X^{(j)} , 1} (b) |; b\in A^{(j+1)} \} = + \infty .
$$
and consequently that
$$
\lim_{j\rightarrow + \infty} \min \{ | \Theta_{X^{(p)} , 1} \Theta_{X^{(p+1)} , 1} \cdots \Theta_{X^{(j)} , 1} (b) |; b\in A^{(j+1)} \} = + \infty .
$$
Therefore we can suppose that $q$, and consequently $u$, is such that each letter of $A^{(p)}$ (the alphabet of $X^{(p)}$) has at least one occurrence in each return word on $u$ of $X^{(p)}$. (We recall that the set of return words on $u$ of $X^{(p)}$ is $\{ \Theta_{X^{(p)},u} (b) = \Theta_{X^{(p)} , 1} \Theta_{X^{(p+1)} , 1} \cdots \Theta_{X^{(q-1)} , 1} (b) ; b\in A^{( q )} \}$.)

The word $\Theta_{X,u^{(p)}} (1)u^{(p)}$ is a prefix of $X$ hence we can choose an integer $l$ such that the word $\Theta_{X,u^{(p)}} (1)u^{(p)}$ is a prefix of $\tau^l (1)$ and $\sigma^l (1)$. Thus the first letter of each image of $\tau_{u^{(p)}}^l$ and $\sigma_{u^{(p)}}^l$ is $1$. 

We set $v=u^{(p)}$ and $\gamma = \tau_{v}^l$. The substitution $(\gamma , A^{(p)} , 1)$ and the prefix $u$ of $X^{(p)}$ (which is the fixed point of $\gamma$) fulfill the hypotheses of Proposition \ref{stepone}. Indeed we chose the integer $l$ to satisfy hypothesis \ref{hypzero}. Hypothesis \ref{hypdeux} is also satisfied because 
$$
\gamma_u = (\tau_{u^{(p)}})_u =  \tau_{u^{(q)}} = \tau_{u^{(p)}} = \gamma ,
$$
where $\gamma_u$ is the return substitution on $u$. Hypothesis \ref{hyptrois} does not set any difficulty. Hypothesis \ref{hypquatre} follows from the choice of $q$.

It is clear that $\sigma^l_{u^{(p)}}$ and $u$ also satisfy the same hypotheses. \hfill $\Box$

\vspace{0,3cm}

\begin{theo}
\label{propvalpr}
If two primitive substitutions have the same non-periodic fixed point, then they have some powers which have the same eigenvalues, except perhaps $0$ and roots of the unity.
\end{theo}
Proof: It follows from Lemma \ref{lemmetech}, Proposition \ref{corotech} and Proposition \ref{relvalpr}.\hfill $\Box$

\subsection{Proof of the main result}

\begin{theo}
\label{maintheo}
Let $X$ be a substitutive sequence arising from  $\tau = (\tau , B , \pl)$, and also from $\sigma = (\sigma ,  C , \pl)$. If $X$ is non-periodic then the dominant eigenvalues of $\tau$ and of $\sigma$ are multiplicatively dependent.
\end{theo}
Proof: Let $A$ be the alphabet of $X$. There exist a morphism $\phi$, from $B$ to $A$, and a morphism $\varphi$ from $C$ to $A$ such that $\phi (\ptf) = \varphi (X_{\sigma}) = X$.

Recall that by Theorem \ref{subretfin}, if a sequence is substitutive then its set of derived sequences is finite. Hence there exist three sequences, $(u^{(n)})_{n\in \pnat}$, $(v^{(n)})_{n\in \pnat}$ and $(w^{(n)})_{n\in \pnat}$, of prefixes of respectively $\ptf$, $X$ and $X_{\sigma}$ such that for all integers $n$ we  have:
\begin{itemize}
\item $\der_{u^{(n)}} (X_{\tau}) = \der_{u^{(n+1)}} (X_{\tau}) $, 
\item $\der_{v^{(n)}} (X) = \der_{v^{(n+1)}} (X) $,
\item $\der_{w^{(n)}} (X_{\sigma}) = \der_{w^{(n+1)}} (X_{\sigma}) $,
\item $\phi(u^{(n)}) = \varphi (w^{(n)}) = v^{(n)}$ and $|v^{(n)}|<|v^{(n+1)}|$.
\end{itemize}
 
Let $n$ be an integer. The images of words by $\phi \Theta_{u^{(n)}}$ are concatenations of return words on $v^{(n)}$. The map $\Theta_{v^{(n)}} : R_{v^{(n)}}^* \rightarrow A^{*}$ being one to one, this allows us to define a morphism $\lambda_n$ by $\Theta_{v^{(n)}} \lambda_n = \phi \Theta_{u^{(n)}}$. In the same way we define the morphism $\gamma_n$ by $\Theta_{v^{(n)}} \gamma_n = \varphi \Theta_{w^{(n)}}$. In the proof of Theorem 21 in \cite{Du}, it is proved that the set $\{ \lambda_n ; n\in \nat \}$, and also the set $\{ \gamma_n ; n\in \nat \}$, are finite. For this reason we can suppose that for all integers $n$ we have $\lambda_n = \lambda_{n+1}$ and $\gamma_n = \gamma_{n+1}$.

Let $i$ be an integer. The sequence $X_{\tau}$ (resp. $X_{\sigma}$) is uniformly recurrent hence, according to Lemma \ref{lemtech}, there exists an integer $j$ larger than $i$ such that each word $wu^{(i)}$, where $w$ is a return word on $u^{(i)}$, has at least one occurrence in each return word on $u^{(j)}$. Consequently we can define a primitive substitution $\delta$ by $\Theta_{u^{(i)}} \delta = \Theta_{u^{(j)}}$. In the same way we define a primitive substitution $\rho $ by $\Theta_{v^{(i)}} \rho = \Theta_{v^{(j)}}$. We have $\rho \lambda_j = \lambda_j \delta$. Indeed
$$
\Theta_{v^{(i)}} \rho \lambda_j = \Theta_{v^{(j)}} \lambda_j = \phi \Theta_{u^{(j)}} = \phi \Theta_{u^{(i)}} \delta = \Theta_{v^{(i)}} \lambda_i \delta = \Theta_{v^{(i)}} \lambda_j \delta.
$$

A standard application of Perron's Theorem (\cite{Ga}, p. 53) shows that $\delta$ and $\rho$ have the same dominant eigenvalue. 

We recall that $\der_{u^{(i)}}(\ptf) = \der_{u^{(j)}}(\ptf)$. Hence $\delta$ has the same fixed point as $\tau_{u^{(i)}}$, that is to say $\der_{u^{(i)}}(\ptf)$. It follows from Theorem \ref{propvalpr} and Proposition \ref{relvalpr} that the dominant eigenvalues of $\delta $ and $\tau $ are multiplicatively dependent.

In the same way we prove that $\rho$ and $\sigma$ have multiplicatively dependent dominant eigenvalues. This completes the proof.\hfill $\Box$

\vspace{0,3cm}

Could we obtain a result analogous to Theorem \ref{propvalpr}? That is to say concerning all eigenvalues. The answer is negative. Here is a counterexample: Let $\tau = (\tau , \{ a,b \} , a)$ and $\sigma = (\sigma , \{ a,b,c \} , a)$ be two substitutions defined respectively by
\[
\left\{
\begin{array}{lll}
a & \rightarrow & abab\\
b & \rightarrow & abbb
\end{array}
\right.
\:\: {\rm and} \:\:
\left\{
\begin{array}{lll}
a & \rightarrow & abab\\
b & \rightarrow & accc\\
c & \rightarrow & abbc
\end{array}
\right.
.
\]

Eigenvalues of the substitution $\tau$ are 1 and 4. Those of $\sigma$ are 1, -2 and 4. Let $\phi : \{ a,b,c \} \rightarrow \{ a,b \}$ be the morphism defined by $\phi (a) = a$ and $\phi (b) = \phi (c) = b$, then $\phi (X_{\sigma}) = \ptf$. The sequence $\ptf$ arises from two substitutions, one has the eigenvalue -2 and the other does not.

\vspace{0,3cm}

To prove the reciprocal of Theorem \ref{maintheo} we need a result  due to D. Lind (Theorem \ref{lind}). A {\it Perron number} is an algebraic integer that strictly dominates all its other algebraic conjugates. It follows easily from Perron's Theorem that the dominant eigenvalue of an integral primitive matrix is a Perron number. The following theorem shows the reciprocal is true.
\begin{theo} 
{\rm (\cite{Li})}
\label{lind}
If $\alpha$ is a Perron number then there exists a primitive integral matrix with dominant eigenvalue $\alpha$.
\end{theo}

\vspace{0,3cm}

Here is the reciprocal of Theorem \ref{maintheo}.

\begin{prop}
Let $Y$ be a periodic sequence on the alphabet $B$ and $\alpha$ a Perron number. There exists an integer $k$ such that $Y$ is $\alpha^k$-substitutive.
\end{prop}
Proof: There exists a word $m$ such that $Y = m^{\omega}$. According to Theorem \ref{lind} there exists an integral primitive matrix $M$ with dominant eigenvalue $\alpha$. There exists an integer $k$ such that:
\begin{enumerate}
\item The matrix $M^k$ has strictly positive coefficients and
\item The sum of the coefficients of any column of $M^k$ is larger than the length of $m$.
\end{enumerate}

It is easy to construct a primitive substitution $\sub$ with associated matrix $M^k$. The dominant eigenvalue of this substitution is $\alpha^k$.

Let $D$ be the alphabet $\{ (b,i); b\in A \:\: , \:\: 0\leq i \leq |m|-1 \}$. We define the morphism $\psi : A \rightarrow D^{+}$ by $\psi(b) = (b,0)\cdots (b,|m|-1)$. The length of an element of $\tau(A)$ is larger than $|m|$. This allows us to define the substitution $\zeta = (\zeta , D , (1,0) )$ in the following way: For all $(b,k)$ of $D$
\begin{center}
\begin{tabular}{llll}
$\zeta((b,k))$   & = & $\psi (\tau (b)_{[k,k]})$           & if $k< |m|-1 $,\\
$\zeta((b,|m|-1))$ & = & $\psi (\tau (b)_{[|m|-1, |\tau(b)|-1]})$ & otherwise.
\end{tabular}
\end{center}

These morphisms are such that $\zeta \psi = \psi \tau$. Hence the substitution $\zeta$ is primitive. Its fixed point is $\psi (X_{\tau})$ and its dominant eigenvalue is $\alpha^k$.

Let $\varphi : D \rightarrow B$ be the letter to letter morphism defined by $\varphi ((b,i)) = m_{[i,i]}$. It is easy to see that $\varphi (X_{\zeta}) = Y$. It follows that $Y$ is $\alpha^k$-substitutive.
\hfill $\Box$

\vspace{0,3cm}

\section{Substitutions sharing the same fixed point}

In this last section we use the circularity of primitive substitutions, proved in \cite{MS,Mo}, to obtain further results about substitutions sharing the same fixed point.
\begin{defi}
Let $\sub$ be a substitution and $x$ a factor of $\ptf$. We say that $(u,w,v)$ is an {\rm interpretation} of $x$ if $x=u\tau(w)v$ and $u$, $v$ are respectively a suffix and a prefix of the image, by $\tau$, of some letters and $w$ is a factor of $\ptf$.
\end{defi}

\begin{defi}
We say that a substitution $\tau$ is {\rm circular} with {\rm synchronization delay D} when: If a factor of $\ptf$ admits two distinct interpretations, $(u,w,v)$ and $(x,y,z)$, and $i$ is an integer such that $|u\tau (w_{[0,i-1]})| >D$ and $| \tau (w_{[i+1,|w|-1]})v| >D$, then there exists an integer $j$ such that $u\tau (w_{[0,i-1]}) = x\tau (y_{[0,j-1]})$ and $w_{i}=y_j$.
\end{defi}

\begin{theo} {\rm (\cite{MS,Mo})}
\label{mignseeb}
A primitive substitution is circular.
\end{theo}

\vspace{0,3cm}

In the following proposition we prove that a primitive substitution is one to one on the set of return words on a sufficiently long prefix of its fixed point.

\begin{prop}
\label{subinj}
Let $\tau$ be a primitive substitution with a non-periodic fixed point $X$. There exists an integer $n_0$ such that for all prefixes $u$ of $X$ of length larger than $n_0$ the substitution $\tau$ is one to one on $L(X) \bigcap \Theta_{X,u} (R_{X,u}^*)$.
\end{prop}
Proof: The substitution $\tau$ is circular with synchronization delay $D$ (Theorem \ref{mignseeb}). According to Lemma \ref{lemtech}, there exists an integer $n_0$ such that for all prefixes $u$ satisfying $|u| > n_0$ the length of all return words on $u$ is larger than $D$.

Let $u$ be a prefix of $X$ larger than ${\rm max} ( n_0 , D)$ and $v$, $w$ be two elements of $L(X) \bigcap \Theta_{X,u} (R_{X,u}^+)$ such that $\tau(v) = \tau(w)$. Let $l$ be the smallest integer $n$ such that $|\tau (v_{[0,n]})| > D$. This integer is smaller than $|u|$ because
$$l\leq |\tau (v_{[0,l-1]})| \leq D < |u| .$$
It follows that $v_{[0,l]}$ is a prefix of $u$. Hence we have $v_{[0,l]} = w_{[0,l]}$. Moreover we have $\tau (vu) = \tau (wu)$ and $|\tau (u)| \geq D$, thus according to Theorem \ref{mignseeb} we obtain $v_{[l+1,|v|-1]} = w_{[l+1,|w|-1]}$ and consequently $v=w$.\hfill $\Box$

\vspace{0,3cm}

\begin{coro}
\label{injsub}
Let $\tau$ be a primitive substitution with a non-periodic fixed point $X$. There exists an integer $n_0$ such that, for all prefixes $u$ of $X$ of length larger than $n_0$, the substitution $\tau_u$ is one to one on $L (X)$.
\end{coro}
Proof: It follows directly from Proposition \ref{subinj}.\hfill $\Box$

\vspace{0,3cm}

To obtain the main result of this section we need an intermediate lemma.
\begin{lemme}
\label{intermed}
Let $\sub$ be a primitive substitution, with fixed point $X$, and $u$ be a prefix of $X$ satisfying the hypothesis of Proposition \ref{stepone}. Let $J$ be an infinite set of positive integers. There exist a subset $I$ of $J$ and a  strictly increasing sequence of positive integers $(l_p)_{p\in I}$ such that for all integers $p$ and $q$, $p<q$, belonging to $I$, we have $\tau^{q-p} = \Theta_u^{l_q - l_p}$.
\end{lemme}
Proof: There exist an infinite subset of $J$, and a strictly increasing sequence $(l_p)_{p\in I}$ such that for all integers $p$ and $q$, $p<q$, we have
$$
\tau^q = \tau^p \Theta_u^{l_q -l_p}.
$$
This follows from equality (\ref{relgamma}) obtained in the proof of Proposition \ref{corotech}. Let $p$ and $q$ be two elements of $I$. It follows from Proposition \ref{derder} that $\ptf$ has a prefix $w$ such that $\Theta_w = \Theta_u^{l_q - l_p}$. Remark that the substitutions $\tau$ and $\tau_w$ are identical. From Lemma \ref{lemtech} and Corollary \ref{injsub} we deduce that we can choose $q$ sufficiently large in order that $\tau_w$ is one to one on $L(X_{\tau_w})$. But $\tau = \tau_w$, hence $\tau$ is one to one on its language. This implies that $\tau^{q-p} = \Theta_u^{l_q - l_p}$.\hfill $\Box$

\vspace{0,3cm}

To end this paper we prove a strong relation between two primitive substitutions sharing the same fixed point.
\begin{prop}
\label{theoprinc1}
Let $\tau$ and $\sigma$ be two primitive substitutions having the same non-periodic fixed point $X$. There exist a prefix $u$ of $X$ and two integers $i$ and $j$ such that
$$
\tau_u^i = \sigma_u^j.
$$
\end{prop}
Proof: There exist a prefix $u$ of $X$ and a prefix $v$ of $\der_u(X)$ such that $\tau_u$ and $v$, and $\sigma_u$ and $v$, both satisfy the hypothesis of Proposition \ref{stepone}. This is Lemma \ref{lemmetech}. It follows from Lemma \ref{intermed} that there exist two integers $i$ and $j$ such that $\tau_u^i = \sigma_u^j$.\hfill $\Box$

\vspace{0,3cm}

With the same hypothesis an equivalent formulation of the previous result is: There exists a prefix $u$ and two integers $i$ and $j$ such that $\tau^i$ and $\sigma^j$ coincide on $L(X) \bigcap \Theta_{X,u} (R^{*}_{X,u})$.

\vspace{0,3cm}

{\bf Aknowledgements} I am grateful to B. Host for the time and the energy he spent on this work. I would like to thank the referees for their valuable advices and constructive comments.

\end{document}